\pdfoutput=1

\documentclass[11pt]{amsart}

\usepackage{epigamath}

\usepackage[english]{babel}

\numberwithin{equation}{section}

\usepackage[shortlabels]{enumitem}
\setlist[enumerate,1]{label={\rm(\arabic*)}, ref={\rm\arabic*}}

\usepackage{amsmath,amsthm, amsfonts, amssymb, mathrsfs, bbm}

\usepackage[capitalize]{cleveref} 

\usepackage{tikz, tikz-cd} 
\usetikzlibrary{positioning}
\usetikzlibrary{decorations.markings}
\makeatletter
\tikzcdset{
	open/.code     = {\tikzcdset{hook, circled};},
	closed/.code   = {\tikzcdset{hook, slashed};},
	open'/.code    = {\tikzcdset{hook', circled};},
	closed'/.code  = {\tikzcdset{hook', slashed};},
	circled/.code  = {\tikzcdset{markwith = {\draw (0,0) circle (.375ex);}};},
	slashed/.code  = {\tikzcdset{markwith = {\draw[-] (-.4ex,-.4ex) -- (.4ex,.4ex);}};},
	markwith/.code ={
		\pgfutil@ifundefined%
		{tikz@library@decorations.markings@loaded}%
		{\pgfutil@packageerror{tikz-cd}{You need to say %
				\string\usetikzlibrary{decorations.markings} to use arrows with markings}{}}{}%
		\pgfkeysalso{/tikz/postaction = {
				/tikz/decorate,
				/tikz/decoration={markings, mark = at position 0.5 with {#1}}}
		}
	},
}
\makeatother

\usepackage{verbatim}
\usepackage[utf8]{inputenc} 
\usepackage[T1]{fontenc} 
\usepackage{comment,braket,xspace,csquotes} 
\usepackage[all,cmtip]{xy} 
\usepackage{multirow}
\usepackage{stmaryrd}

\usepackage{tikz, tikz-cd} 
\usetikzlibrary{positioning}
\usetikzlibrary{decorations.markings}
\usepackage{xcolor}

\theoremstyle{plain}
\newtheorem{proposition}{Proposition}[section]
\newtheorem{lemma}[proposition]{Lemma}
\newtheorem{theorem}[proposition]{Theorem}
\newtheorem*{theorem*}{Theorem}
\newtheorem{corollary}[proposition]{Corollary}
\newtheorem*{corollary*}{Corollary}

\theoremstyle{definition}
\newtheorem{definition}[proposition]{Definition}

\theoremstyle{remark}
\newtheorem{remark}[proposition]{Remark}

\usetikzlibrary{decorations.markings}

\makeatletter
\tikzcdset{
	open/.code     = {\tikzcdset{hook, circled};},
	closed/.code   = {\tikzcdset{hook, slashed};},
	open'/.code    = {\tikzcdset{hook', circled};},
	closed'/.code  = {\tikzcdset{hook', slashed};},
	circled/.code  = {\tikzcdset{markwith = {\draw (0,0) circle (.375ex);}};},
	slashed/.code  = {\tikzcdset{markwith = {\draw[-] (-.0ex,-.4ex) -- (.0ex,.4ex);}};},
	markwith/.code ={
		\pgfutil@ifundefined%
		{tikz@library@decorations.markings@loaded}%
		{\pgfutil@packageerror{tikz-cd}{You need to say %
				\string\usetikzlibrary{decorations.markings} to use arrows with markings}{}}{}%
		\pgfkeysalso{/tikz/postaction = {
				/tikz/decorate,
				/tikz/decoration={markings, mark = at position 0.5 with {#1}}}
		}
	},
}
\makeatother

\newcommand{\rH}{\mathrm H}

\newcommand{\cL}{\mathcal L}
\newcommand{\cM}{\mathcal M}

\newcommand{\cO}{\mathcal O}

\newcommand{\bZ}{\mathbb{Z}}

\newcommand{\bA}{\mathbb{A}}

\newcommand{\bQ}{\mathbb{Q}}
\newcommand{\ubQ}{\underline\bQ}

\newcommand{\bC}{\mathbb{C}}

\newcommand{\bR}{\mathbb{R}}
\newcommand{\bS}{\mathbb{S}}

\newcommand{\bG}{\mathbb{G}}

\makeatletter
\newcommand{\oset}[3][0.2ex]{%
	\mathrel{\mathop{#3}\limits^{
			\vbox to#1{\kern-2\ex@
				\hbox{$\scriptstyle#2$}\vss}}}}
\makeatother

\newcommand{\pt}{\mathsf{pt}}

\newcommand{\inv}{^{-1}}

\usetikzlibrary{decorations.markings} 
\tikzset{
  closed/.style = {decoration = {markings, mark = at position 0.5 with { \node[transform shape, xscale = .8, yscale=.4] {/}; } }, postaction = {decorate} },
  open/.style = {decoration = {markings, mark = at position 0.5 with { \node[transform shape, scale = .7] {$\circ$}; } }, postaction = {decorate} }
}

\DeclareMathOperator{\Aut}{Aut}
\DeclareMathOperator{\Hom}{Hom}

\DeclareMathOperator{\Pic}{Pic}
\DeclareMathOperator{\CH}{CH}

\DeclareMathOperator{\Spec}{Spec}

\DeclareMathOperator{\Ima}{Im}

\makeatletter
\let\save@mathaccent\mathaccent
\newcommand*\if@single[3]{%
	\setbox0\hbox{${\mathaccent"0362{#1}}^H$}%
	\setbox2\hbox{${\mathaccent"0362{\kern0pt#1}}^H$}%
	\ifdim\ht0=\ht2 #3\else #2\fi
}
\newcommand*\rel@kern[1]{\kern#1\dimexpr\macc@kerna}
\newcommand*\widebar[1]{\@ifnextchar^{{\wide@bar{#1}{0}}}{\wide@bar{#1}{1}}}
\newcommand*\wide@bar[2]{\if@single{#1}{\wide@bar@{#1}{#2}{1}}{\wide@bar@{#1}{#2}{2}}}
\newcommand*\wide@bar@[3]{%
	\begingroup
	\def\mathaccent##1##2{%
		\let\mathaccent\save@mathaccent
		\if#32 \let\macc@nucleus\first@char \fi
		\setbox\z@\hbox{$\macc@style{\macc@nucleus}_{}$}%
		\setbox\tw@\hbox{$\macc@style{\macc@nucleus}{}_{}$}%
		\dimen@\wd\tw@
		\advance\dimen@-\wd\z@
		\divide\dimen@ 3
		\@tempdima\wd\tw@
		\advance\@tempdima-\scriptspace
		\divide\@tempdima 10
		\advance\dimen@-\@tempdima
		\ifdim\dimen@>\z@ \dimen@0pt\fi
		\rel@kern{0.6}\kern-\dimen@
		\if#31
		\overline{\rel@kern{-0.6}\kern\dimen@\macc@nucleus\rel@kern{0.4}\kern\dimen@}%
		\advance\dimen@0.4\dimexpr\macc@kerna
		\let\final@kern#2%
		\ifdim\dimen@<\z@ \let\final@kern1\fi
		\if\final@kern1 \kern-\dimen@\fi
		\else
		\overline{\rel@kern{-0.6}\kern\dimen@#1}%
		\fi
	}%
	\macc@depth\@ne
	\let\math@bgroup\@empty \let\math@egroup\macc@set@skewchar
	zen@everymath{\mathgroup\macc@group\relax}%
	\macc@set@skewchar\relax
	\let\mathaccentV\macc@nested@a
	\if#31
	\macc@nested@a\relax111{#1}%
	\else
	\def\gobble@till@marker##1\endmarker{}%
	\futurelet\first@char\gobble@till@marker#1\endmarker
	\ifcat\noexpand\first@char A\else
	\def\first@char{}%
	\fi
	\macc@nested@a\relax111{\first@char}%
	\fi
	\endgroup
}
\makeatother

\newcommand{\lra}{\longrightarrow}

\EpigaVolumeYear{8}{2024} \EpigaArticleNr{11} \ReceivedOn{September 28, 2023}
\InFinalFormOn{January 26, 2024}
\AcceptedOn{February 25, 2024}

\title{Ng\^o's support theorem and polarizability  of quasi-projective commutative group schemes}
\titlemark{Support theorem and polarization}

\author{Giuseppe Ancona}
\address{IRMA, Strasbourg, France}
\email{ancona@math.unistra.fr}
\author{Drago\c s Fr\u a\c til\u a}
\address{IRMA, Strasbourg, France}
\email{fratila@math.unistra.fr}

\authormark{G.~Ancona and D.~Fr\u{a}\c{t}il\u{a}}

\AbstractInEnglish{In this short note, we prove that any commutative group scheme $G\to B$ over an arbitrary base scheme $B$ of finite type over a field   with connected fibers and admitting a relatively ample line bundle is polarizable in the sense of Ng\^o.
More precisely, we associate a polarization to any relatively ample
line bundle on $G\to B$.  This extends the applicability of Ng\^o's
support theorem to new cases, such as Lagrangian fibrations with
integral fibers, and has consequences for the construction of algebraic
classes.}



\acknowledgement{This research was partly supported by  grant ANR--18--CE40--0017 of the Agence Nationale de la Recherche.}

\begin{document}


\maketitle

\begin{prelims}

\DisplayAbstractInEnglish

\bigskip

\DisplayKeyWords

\medskip

\DisplayMSCclass

\end{prelims}


\newpage

\setcounter{tocdepth}{1}

\tableofcontents


\section{Introduction}

Let $\pi\colon G\to B$ be a commutative group scheme of relative dimension $d$ with connected fibers over a scheme $B$ of finite type over $\bC$. 
The relative Tate module is 
\begin{align*}\label{eq def tate module}
T(G):=R^{2d-1}\pi_!\ubQ_G(d),
\end{align*} 
where we denote by $(d)$ the Tate twist and $\ubQ_G$ is the constant sheaf on $G$.

By  Chevalley's structure theorem, for all $b \in B$, the fiber  $G_b$ of $\pi$ over $b$ sits in an exact sequence
\begin{equation}\label{chevalley}
	1 \longrightarrow L_b \longrightarrow G_b \longrightarrow A_b  \longrightarrow   1,
\end{equation}
 where $L_b$ is the unique connected affine   subgroup of $G_b$ such that the quotient group~$A_b$ is an abelian variety.

The following is the definition of polarization after \cite[Section~7.1.4]{ngo2010lemme}. 

\begin{definition}\label{d polariz}
	A polarization on $G$ is an alternating bilinear map
	\[ \eta\colon T(G)\otimes T(G)\longrightarrow \ubQ_B(1) \]
	such that for every geometric point $b\in B$, the induced bilinear form
\begin{equation}\label{eq polariz}
	\eta_b\colon T(G_b)\otimes T(G_b)\longrightarrow \bQ(1)
\end{equation}
	has kernel exactly $T(L_b)$, where $L_b$ is as in \eqref{chevalley}.
	
	If such a polarization exists, we say that $T(G)$, or simply $G$, is polarizable.
\end{definition}

Our main result is the following.

\begin{theorem}\label{t main th}
	Any quasi-projective group scheme is polarizable. 
	More precisely, 
	let $\cL$ be  a relatively ample line bundle  on $G\to B$.
	Then its first Chern class   induces a polarization
	\[ \eta_{\cL}\colon T(G)\otimes T(G)\longrightarrow \ubQ_B(1).\]
\end{theorem}

The original definition of Ng\^o uses $\ell$-adic coefficients instead of rational coefficients and makes sense over very general bases. Our result also holds  in this setting (see \cref{S l-adic}).

Recently, Ng\^o's support theorem has been extensively used to study Lagrangian fibrations. One of its crucial hypotheses is the polarizability of the underlying group scheme. 
Hence, our \cref{t main th} allows one to apply Ng\^o's support theorem to more general Lagrangian fibrations, and in particular we can deduce the following result. 

\begin{corollary}\label{cor HK}
Let $X$ be a projective hyper-k\"ahler variety and $f\colon X \rightarrow B$ be a Lagrangian fibration with integral fibers. Then all perverse sheaves appearing in the decomposition theorem for $f$ have dense support.
\end{corollary}

\begin{remark}
During the final stages of this work, we  learned that similar results with different proofs have been obtained by  Mark de Cataldo, Roberto Fringuelli, Andr\'es Fernandez-Herrero and Mirko Mauri. Their article is in preparation. We warmly thank them for the friendly exchanges during our parallel work.
\end{remark}

\subsection*{Idea of the proof of \cref{t main th}} The construction of $\eta_\cL$ is purely formal from the adjunction $(R\pi_!, \pi^!)$. In particular, base change ensures functoriality by pullback on the base, and the problem is therefore reduced to the absolute case $B=\pt=\{b\}$ and $\cL_b$ an ample line bundle on the commutative algebraic group $G_b$. 
Such a reduction step is the main reason why we construct a polarization using Chern classes instead of constructing it using  purely Hodge-theoretical techniques. Notice that it is easy to show that every $G_b$ is polarizable using Hodge theory; however, it is unclear in general how to put those polarizations in a family. 
The use of Chern classes allows us to give a global definition of a pairing which restricts well on each fiber. The point then becomes to check that the restriction on each fiber is a polarization.

Let  us now write $G=G_b$ and  consider Chevalley's structure theorem applied to $G$: 
\[1\lra L\lra G\stackrel{p}{\lra} A\lra 1.\] 
One shows that $p$ induces a surjective map
\[p^*\colon \Pic(A)\lra \Pic(G).\]
We moreover show, and it is essential, that if $\cL:=p^*\cM$ is ample on $G$, then $\cM$ is ample as well (see \cref{p L e M ample}, or \cref{p ample on G non-deg H on A} for a weaker version). This reduces the problem to the case of an abelian variety, where one has to show that the pairing is non-degenerate.

When  $G=A$, we check that the pairing constructed in the paper is  related to  more classical polarizations from Hodge theory, which are known to be non-degenerate. This is possible through the explicit analytic description of   Chern classes on abelian varieties given by the
 Appell--Humbert theorem (see \cite[Section~I.2, p.~20]{MumfordTata}).
 An alternative algebraic argument is provided in \cref{S l-adic}.

\subsection*{Organization of the paper}	
In Section \ref{S constr polariz}, we explain the construction of the polarization associated with a line bundle. The special absolute case of complex abelian varieties is studied in Section \ref{abel var}. The general absolute case is explained in Section \ref{absolute case}. The reduction to the absolute case is explained in Section \ref{end proof}, where one will also find the proofs of \cref{t main th} and \cref{cor HK}. 
In \cref{S l-adic}, we explain how to modify the arguments in order to extend the results from complex bases to general schemes.

\subsection*{Acknowledgments}	
We thank Marian Aprodu, Olivier Benoist, Michel Brion, Mattia Cavicchi, Mark de Cataldo and Giulia Sacc\`a for   their interest in our work and for  fruitful discussions. We warmly thank the   referee for a careful reading.

\section{Construction of the polarization}\label{S constr polariz}
Let $\pi \colon G\to B$ be a commutative group scheme of relative dimension $d$. Recall that we defined the Tate module to be $T(G):=R^{2d-1}\pi_!\ubQ_G(d)$.
In this section, we construct, for any class in cohomology $\omega\in \rH^2(G,\bQ)(1)$, a map
\[ \eta_\omega\colon \Lambda^2 T(G)\lra \ubQ_B(1).\]
We will apply it in later sections to $\omega = c_1(\cL)$ for $\cL$ a relatively ample line bundle on $G\to B$.

 We can view $\omega$  as a map in the derived category of constructible sheaves on $G $
\[ \omega\colon \ubQ_{ G}\lra \ubQ_{G }[2](1).\]

For a smooth map of finite type schemes $f\colon X\to Y$ of relative dimension $d$, we have $f^* = f^![-2d](-d)$. This is a relative version of Poincar\'e duality (see \cite[Theorem XVIII.3.2.5]{SGA4t3} and its proof).
Applying this to the smooth map $\pi $, we get
\[\ubQ_{G} = \pi ^!\ubQ_B[-2d](-d),\]
and $\omega$ can now be written  as
\[ \omega\colon \ubQ_{G }\lra \pi^!\ubQ_B[-2d+2](-d+1).\]
By the adjunction $(R\pi_!, \pi^!)$, we can view $\omega$ equivalently as
 \begin{align}\label{eq omega}
	 \omega\colon R\pi_!\ubQ_{G }\lra  \ubQ_B[-2d+2](-d+1),
\end{align}
and taking the $2d-2$ cohomology sheaf in Equation \eqref{eq omega} gives the map 
\begin{align}\label{eq constr eta}
	\eta_\omega\colon R^{2d-2}\pi_!\ubQ_G  \lra \ubQ_B(-2d+1).
\end{align}
Through the canonical identification $R^{2d-2}\pi_!\ubQ_G=  \Lambda^2 R^{2d-1}\pi_!\ubQ_G$  (see for example \cite{AHP}), we have indeed constructed
\begin{align}\label{eq final eta}
	 \eta_\omega\colon \Lambda^2 T(G)\lra \bQ_B(1).
\end{align}
By its very definition, this construction is compatible with changing the base $B$. 

\begin{lemma}\label{l pullback polariz}
	Let $f\colon B'\to B$ be a map from a finite type scheme $B'$, and denote by $\pi'\colon G'\to B'$ the base change to $B'$. Let $\omega\in H^2(G,\bC)(1)$, and let $\omega'\in H^2(G',\bC)(1)$ be its pullback.
	Then we have $f^*(\eta_\omega) = \eta_{\omega'}$; i.e., the following diagram commutes: 
	 \[\begin{tikzcd}
		f^*T(G)\otimes f^*T(G)\ar[r,"f^*(\eta_\omega)"]\ar[d,"\sim"] & f^*\ubQ_B(1)\ar[d]\\
		T(G')\otimes T(G')\ar[r,"\eta_{\omega'}"] & \ubQ_{B'}(1)\rlap{,}
		\end{tikzcd} \]
	where the left vertical isomorphism comes from base change.
\end{lemma}
 
\section{The case of abelian varieties: The Appell--Humbert theorem}\label{abel var}
The goal of this section is to prove that $\eta_\omega$ as in  \eqref{eq final eta} is  non-degenerate when $\omega$ is the Chern class of   an ample line bundle $\cL$ on an abelian variety $A$.
After some preliminaries, we will recall (following \cite[Section~I.2]{MumfordTata})
the main ingredient: the Appell--Humbert theorem.

Let $A$ be a complex abelian variety,  and write it as $V/\Gamma$, where $V=T_eA$ and $\Gamma = \rH_1(A,\bZ)$. The quotient is taken in the realm of complex analytic varieties.

Any line bundle $\cL$ on $A$ becomes trivial when pulled back to $V$. 
To recover $\cL$, some descent datum must be specified; in this situation, it can be made very explicit in terms of linear algebra on $V$ and $\Gamma$. 
Consider pairs
$(H,\rho)$, where
\[H\colon V\otimes_\bR V\lra \bC\] 
is a hermitian form and 
\[\rho\colon \Gamma\lra \bS^1\]
is an $H$-pseudo-character; \textit{i.e.}, $\rho$ satisfies
\begin{align}\label{eq rho H pseudo char}
	 \rho(u+v) = e^{i\pi E(u,v)}\rho(u)\rho(v)\quad \text{for all }u,v\in\Gamma, 
\end{align}
where $E:=\Ima H\colon V\otimes_\bR V\to \bR$ is the imaginary part of $H$.

With such a datum $(H,\rho)$, assuming that $E(\Gamma\times\Gamma)\subset \bZ$, we can associate the line bundle $\cL(H,\rho)$ on $A$ given by
\begin{align}\label{eq constr of L(H,rho)} 
	\cL(H,\rho):=( V\times \bC)/\Gamma, 
	\end{align}
where the action of $\Gamma$ is given by
\[ u\cdot (v,z):=(u+v,\rho(u)e^{\pi H(v,u)+\frac12\pi H(u,u)}z).\]

We list some useful properties whose proofs are immediate. 

\begin{lemma} With the above notation, the following hold:
\begin{enumerate}
	\item If\, $(H_1,\rho_1)$ and $(H_2,\rho_2)$ are two pairs as above, then
\[ \cL(H_1,\rho_1)\otimes \cL(H_2,\rho_2) = \cL(H_1+H_2,\rho_1\rho_2). \]
\item If $f\colon A'\to A$ is a morphism of abelian varieties with $A'=V'/\Gamma'$, then
\[ f^*(\cL(H,\rho)) = \cL(f^*H,f^*\rho), \]
where $f^*H$ is the hermitian form on $V'$ induced from $H$ through the linear map $V'\to V$ and similarly for $f^*\rho$. 
\end{enumerate}
\end{lemma}

For the comfort of the reader, let us also recall the following relationship between $H$ and $E$.

\begin{lemma}[\textit{cf.} \protect{\cite[Section~I.2, p.~20]{MumfordTata}}]\label{l hermitian vs skew-sym}
	There is a bijection between hermitian forms $H$ on $V$ and real skew-symmetric forms $E$ on $V$ satisfying the identity $E(ix,iy) = E(x,y)$, which is given by 
	\begin{align}
		E(x,y) &= \Ima H(x,y),\\
		H(x,y)& = E(ix,y)+iE(x,y).
	\end{align}
Moreover, the kernel of $H$ is equal to the kernel of $E$. In particular, $H$ is non-degenerate if and only if $E$ is non-degenerate.
\end{lemma}

The following result describes the N\'eron--Severi group 
\[ NS(A) := \textrm{Image}(\rH^1(A,\cO_A^\times)\lra \rH^2(A,\bZ)) \]
under the identification
\[ \rH^2(A,\bZ) = (\Lambda^2\Gamma)^\vee.\]

\begin{theorem}[Appell--Humbert, \textit{cf.} \protect{\cite[Section~I.2, p.~20]{MumfordTata}}]\label{t AH}
Any line bundle on $A$ is of the form $\cL(H,\rho)$ for a unique pair $(H,\rho)$ as above. Moreover, we have isomorphic short exact sequences
\[
\begin{tikzcd}
	0\ar[r] & \Hom(\Gamma,\bS^1)\ar[r]\ar[d]&
	\left\{\begin{array}{c}(H,\rho)\\  \text{satisfying }\\\eqref{eq rho H pseudo char}\end{array} \right\}\ar[d]\ar[r] & \left\{
	\begin{array}{cc}H \text{ hermitian}\\ \text{on $V$ s.t. } \\
		\Ima(H)(\Gamma\times\Gamma)\subseteq\bZ
	\end{array}\right\} \ar[d,"\gamma"]\ar[r]& 0\\
	0\ar[r] & \Pic^0(A)\ar[r] & \Pic(A) \ar[r,"c_1(-)"] & NS(A) \ar[r] &0\rlap{,}
\end{tikzcd}
\]
where $\gamma(H) = E:=\Ima(H)$ is the imaginary part of $H$.
\end{theorem}

\begin{corollary}\label{c chern class is E}
With the above notation, we have
\[c_1(\cL(H,\rho)) = E\in (\Lambda^2\Gamma)^\vee=  \rH^2(A,\bZ), \quad \text{where }E=\Ima(H).\]
\end{corollary}

\begin{theorem}[Lefschetz's theorem, \textit{cf.} \protect{\cite[Section~I.3, p.~28]{MumfordTata}}]\label{lefsc}
	Keeping the above notation, we have that $\cL(H,\rho)$ is ample if and only if $H$ is positive definite.
\end{theorem}

We now have the ingredients to prove the main result of the section, which is a special case of \cref{t main th}.

\begin{theorem}\label{t polariz for abelian varieties}
Let $\cL(H, \rho)$ be an ample line bundle on the abelian variety $A=V/\Gamma$.
Then the map $\eta_\omega$ constructed in \eqref{eq final eta} associated with $\omega = c_1( \cL(\rho,H))$, 
\[ \eta_\omega\colon \Lambda^2 T(A)\lra \bQ(1), \]
is non-degenerate.
\end{theorem}

\begin{proof}
	Recall that, starting with a class $\omega \in H^2(A)$, the construction of the pairing $\eta_\omega$ is based on the identification $H^2(A)=\Lambda^2 H^1(A)$ and on Poincar\'e duality.
	 On the other hand, notice that the description of $c_1(\cL(H,\rho))$ in \cref{c chern class is E} is also based on the equality $H^2(A)=\Lambda^2 H^1(A)$ and on duality.
	
	Now we have the identification $T(A)= \Gamma$, as they are both the dual of $H^1(A)$. Under this identification, the pairing $\eta_{ c_1( \cL(\rho,H))}$ coincides with $E= \Ima(H)$, by \cref{c chern class is E}. As $\cL(H, \rho)$ is ample, the hermitian form $H$ is non-degenerate by \cref{lefsc}, which implies that $E$ is non-degenerate as well by \cref{l hermitian vs skew-sym}.
\end{proof}

\section{Picard groups of commutative algebraic groups}\label{absolute case}
Let $G$ be a connected, commutative group scheme of finite type over $\bC$, and use Chevalley's structure theorem to get
\[ 1\lra L\lra G\stackrel{p}{\lra} A\lra 1 \]
with $L$ an affine algebraic group and $A$ an abelian variety.
The purpose of this section is to relate the Picard group of $G$ with that of $A$. This is achieved in \cref{p ample on G non-deg H on A},  which is the only result of this section using analytic methods. In \cref{p L e M ample}, we propose an algebraic alternative for it.

\begin{proposition}\label{p surj on Picard}
	With the above notation, the map $p\colon G\to A$ induces a surjection
	\[ p^*\colon \Pic(A)\lra\Pic(G).\]
\end{proposition}

\begin{proof} 
If $L$ has a non-trivial unipotent radical, then by $\bA^1$-homotopy, 
we clearly have $\Pic(G)\simeq \Pic(G/L_u)$, where $L_u$ is the
unipotent radical of $L$.  We can therefore assume that $L$ has no
unipotent subgroups, in other words, that $L\simeq \bG_m^r$ for some
$r\ge 1$.
	
	Consider the natural action of $\bG_m^r$ on $\bA^r$, and take the associated $\bA^r$-bundle
	\[ \overline G:=G\times^{L}\bA^r := (G\times \bA^r)/L\]
	giving the following commutative diagram, with $j$ an open immersion: 
	\[ \begin{tikzcd}
		G \ar[r,hook,"j",circled]\ar[d,"q",swap] & \overline G\ar[dl,"\overline q"] \\
		A\rlap{.}
	\end{tikzcd} \]
	Consider the restriction  map
	\[ j^*\colon \CH^1(\overline G)\lra \CH^1(G). \]
        It is surjective since any cycle on $G$ is the restriction of its closure in $\overline G$.
	Moreover, since $\overline q^*\colon \CH^1(A)\to \CH^1(\overline G)$ is an isomorphism by $\bA^1$-homotopy, we deduce that $q^*\colon \CH^1(A)\to \CH^1(G)$ is also surjective.
\end{proof}

For the next result, we need to use some sort of numerical criterion for ample divisors on open varieties.\footnote{We warmly thank Olivier Benoist for telling us about these results.}

\begin{proposition}\label{p ample on open subvar}
	Let $X$ be a smooth projective variety, and let $U\subset X$ be an open subvariety.
	Let $\cL_1,\cL_2$ be two line bundles on $X$ such that their first Chern classes coincide $c_1(\cL_1)=c_1(\cL_2)$.
	If $\cL_1|_U$ is ample, then $\cL_2|_U$ is also ample.
\end{proposition}

\begin{proof}
	It is basically \cite[Proposition 7]{OB_qproj}.
	In \textit{loc.~cit.}, the statement is for algebraic equivalence. However, for divisors on smooth projective varieties, algebraic equivalence coincides with  numerical equivalence by a theorem of Matsusaka (1956); see \cite[Corollary 1]{kleiman1966toward}.
\end{proof}

In our context, we will apply it as follows

\begin{lemma}\label{l chern=0 quas-aff}
	Let $X$ be a smooth quasi-projective variety and $\cL$ an ample line bundle on $X$.
	Suppose the first Chern class of $\cL$ vanishes.
	Then $X$ is quasi-affine.
\end{lemma}

\begin{proof}
	Let $X\subset \overline X$ be a smooth compactification with complement a normal crossing divisor (this is possible in characteristic zero, by Nagata's theorem and the resolution of singularities).
	Denote by $Z$ the divisor at infinity.
	
	Extend the line bundle $\cL$ to a line bundle $\overline\cL$ over $\overline X$ by taking the closure of the corresponding Weil divisor.
	The Chern class of $\overline\cL$ is going to be supported on $Z$ because $c_1(\overline\cL|_U)=c_1(
	\cL)=0$.
	Hence, there is a divisor supported on $Z$, call it $D$, such that $c_1(\overline \cL) = c_1(\cO_{\overline X}(D))$.
	
	We can now apply \cref{p ample on open subvar} and deduce that $\cO_X=\cO_{\overline X}(D)|_X$ is ample on $X$; this means precisely that $X$ is quasi-affine (see \cite[\href{https://stacks.math.columbia.edu/tag/01QE}{Tag 01QE}]{stacks-project}).
\end{proof}

The next lemma is contained in \cite[Proposition VIB.11.11]{SGA3}, but we include a short proof for completeness.

\begin{lemma}\label{l qaffine group is affine}
	Let $G$ be a commutative group scheme. If\, $G$ is quasi-affine, then $G$ is affine.
\end{lemma}

\begin{proof}
	Let us put $G':=\Spec(\cO(G))$. We have a morphism of group schemes
	$G\to G'$ that is moreover an open dense embedding because of the assumption that $G$ is quasi-affine.
	Since images of group morphisms are closed (see for example \cite[\href{https://stacks.math.columbia.edu/tag/047T}{Lemma 047T}]{stacks-project}), we deduce that $G=G'$.
\end{proof}

 Recall the notation $\cL(H,\rho)$ from  \cref{t AH} and the preceding paragraphs.

\begin{proposition}\label{p ample on G non-deg H on A} 
	  Assume that the  line bundle $p^*\cL(H,\rho)$ on $G$ is ample.    Then $H$ is a non-degenerate hermitian form.
\end{proposition}

\begin{proof}
Let us use the complex uniformization of $A$: put $V=T_eA$ and $\Gamma:=\rH_1(A,\bZ)$.
We have $V/\Gamma \simeq A$ (quotient in the realm of complex analytic varieties).

Put $N:=\ker(H)$. It is a complex vector subspace of $V$ which moreover has the property that $N\cap\Gamma$ is cocompact in $N$, as the imaginary part $E:=\Ima(H)$ of $H$ takes integral values on $\Gamma$ (see the condition on $\Ima(H)$ in \cref{t AH}).

Therefore, we can consider $A':=N/(N\cap \Gamma)$ as a complex subtorus of $A$. By Chow's lemma, we have that $A'$ is projective, hence an abelian variety.
Moreover, the line bundle $\cL(H,\rho)$ restricted to $A'$ is precisely $\cL(0,\rho|_{\Lambda\cap N})$, again by \cref{t AH}.
In particular, using \cref{c chern class is E}, we have that its Chern class is zero.

The restriction of $p\colon G\to A$ to $A'$ gives a commutative group scheme with Chevalley decomposition
\[ 1\lra L\lra G'\lra A'\lra 1 \]
with $G':=p\inv(A')$ and such that $(p^*\cL(H,\rho))|_{G'}$ is ample.
The Chern class of this line bundle is zero by functoriality; hence $G'$ is quasi-affine by \cref{l chern=0 quas-aff} and therefore affine by \cref{l qaffine group is affine}.
It follows that $A'=1$ and therefore that $N=\{0\}$ or, in other words, that $H$ is non-degenerate.
\end{proof}

\begin{remark}
As already  pointed out,  \cref{p ample on G non-deg H on A}    is the only result in  this section using analytic methods. In \cref{p L e M ample}, we propose an algebraic alternative for it.
\end{remark}

\section{Proof of the main theorem}\label{end proof}
We are now ready to give the proof of our main result \cref{t main th} and of its application \cref{cor HK}. For the convenience of the reader, we recall the statements.

\begin{theorem}\label{main thm}
Let $B$ be a scheme of finite type over the complex numbers. Let 
	 $\pi\colon G\to B$ be a commutative group schemes with  connected fibers. Assume that there is a    relatively ample line bundle   $\cL$ on $\pi \colon G\to B$.    
	Let $\omega$ be the first Chern class of  $  \cL.$ Then 
	\[ \eta_{\omega}\colon T(G)\otimes T(G)\lra \ubQ_B(1)\]
	\textup{(}as constructed in Section \ref{S constr polariz}\,\textup{)} is a polarization in the sense of \cref{d polariz}.
\end{theorem}

\begin{proof}
	By the functoriality of the construction, see \cref{l pullback polariz}, we are reduced to the case $B=\Spec(\bC)$.	
	
	Let $1\to L\to G\to A\to 1$ be the Chevalley decomposition of $G$ with $L$ an affine algebraic group and $A$ an abelian variety. We want to check that $\eta_{\omega}$ has kernel precisely $T(L)$. Notice that for weight reasons, the kernel contains $T(L)$; hence we have to show that the induced pairing on $T(A)$ is non-degenerate.
		
	By \cref{p surj on Picard}, $\cL$ is the pullback of a line bundle on $A$. Hence, by \cref{t AH}, we can write $\cL=p^*\cL(H,\rho)$, where $p\colon G\to A$ is the projection. By  \cref{l pullback polariz}, the pairing we want to study on $T(A)$ is the one induced by $ \cL(H,\rho)$. This is nothing but $\Ima H$, by \cref{c chern class is E}. On the other hand,  \cref{p ample on G non-deg H on A} shows that $H$ is non-degenerate, hence so is $\Ima H$ by virtue of
	\cref{l hermitian vs skew-sym}.
\end{proof}

\begin{corollary}\label{lagrangian}
Let $X$ be a projective hyper-k\"ahler variety and $f\colon X \rightarrow B$ be a Lagrangian fibration with integral fibers. Then all perverse sheaves appearing in the decomposition theorem for $f$ have dense support.
\end{corollary}

\begin{proof}
In \cite{AF16}, under the same hypotheses, the authors construct a group scheme $\pi\colon G \to B$ as a subgroup of $\Aut(f)$ and verify that all  hypotheses needed to apply Ng\^o's support theorem are satisfied, except possibly the polarization (see also \cite[Proposition 9.5]{ACLS}). On the other hand, under the hypothesis of irreducible fibers,  Ng\^o's support theorem implies that the supports of the perverse sheaves are dense. We are then reduced to showing that $\pi$ is polarizable; hence, by \cref{main thm}, we have to show that $\pi$ admits a relatively ample line bundle. 

By the construction of $G$, it is enough to show that $\Aut(f)$ admits a relatively ample line bundle.
If we associate with each automorphism its graph, we obtain an embedding  of $\Aut(f)$ in the Hilbert scheme of $X \times_B X$. The latter is a disjoint union of projective schemes over $B$, so we can conclude.
\end{proof}

\begin{remark}\label{ACLS}
This corollary has applications to the construction of algebraic classes on Lagrangian fibrations; see the main results of \cite{ACLS}.
\end{remark}

\section{The $\ell$-adic setting}\label{S l-adic}
In this section, we extend our main theorem to group schemes over more general bases, in particular in positive characteristic and in mixed characteristic.

Let $\ell$ be a prime number and $B$ be a base scheme such that   $\ell$ is  invertible on $B$. This assumption is needed in order to have $\ell$-adic cohomology. 
Let $\pi\colon G\to B$ be a commutative group scheme with  connected fibers of relative dimension $d$.
The Tate module is defined as before as 
\[ T  (G):=R^{2d-1}\pi_!\ubQ_\ell (d). \]

\begin{definition} (\textit{cf.} \cite[Section 7.1.4]{ngo2010lemme})
	A polarization on $\pi\colon G\to B$ is a pairing on the Tate module
	\[ \eta\colon \Lambda^2 T(G)\lra \ubQ_\ell(1)\] 
	such that for every geometric point $b$ of $B$, if
	\[ 1\lra L_b\lra G_b\lra A_b\lra 1 \]
	is the Chevalley decomposition, then the pairing $\eta_b$ induces a non-degenerate pairing
	\[ \eta_b\colon T(A_b)\otimes T(A_b)\lra \bQ_\ell(1).\]
	
	If such a polarization exists, we say that $G$ is polarizable.
\end{definition}

\begin{theorem}\label{t:main general}
	Assume that $B$ is quasi-compact and quasi-separated and that $\pi$ is quasi-projective. Then $G$ is polarizable. More precisely, 
	   a relative ample line bundle on $\pi $ induces a polarization through its  first Chern class.
\end{theorem}

\begin{remark}
	The assumptions of quasi-compact and quasi-separated are those in \cite[Chapter~XVIII]{SGA4t3}   that are needed for duality. Ng\^o's definition is, anyway, in the context of $B$ being of finite type over a field. 
\end{remark}

\begin{proof}	
	The proof of \cref{t main th} applies to the $\ell$-adic context \textit{mutatis mutandis}, up to three exceptions.
	
	First,  the construction of the pairing $\eta_\omega$ associated with a cohomology class $\omega$ used the duality $\pi^*=\pi^![-2d](-d)$. In characteristic zero, such a duality holds as $\pi$ is always smooth, but in our generality, we need an extra argument. 
	
	\smallskip
	
	Notice that $\pi$ is compactifiable in the sense of \cite[D\'efinition XVII.3.2.1]{SGA4t3}, as it is quasi-projective. This gives a trace map relating $\pi^*$ and $\pi^![-2d](-d)$ by \cite[Equation XVIII.(3.2.1.2)]{SGA4t3}. 
	As the trace map is compatible with base change (by \textit{loc.~cit.}), it is enough to check that it is an isomorphism on geometric fibers. 
	The fibers are not smooth, but they become smooth with the reduced scheme structure (as they are algebraic groups). 
	Since $\ell$-adic cohomology is insensitive to the scheme structure, we can conclude that the trace map $\pi^*\to \pi^![-2d](-d)$ is an isomorphism.
	
	\smallskip	
	
	Second, the reduction to the case of an abelian variety used analytic methods (\cref{p ample on G non-deg H on A}). It can be made purely algebraic by using the functoriality of the pairing \cref{l pullback polariz} and \cref{p L e M ample} below.
	
		\smallskip	
		
	Third, 
  the proof in the case of an abelian variety (see \cref{t polariz for abelian varieties}) was also of analytic nature, relying on the Appell--Humbert theorem and analytic uniformization. 
	We give here an alternative argument of algebraic nature, which is less precise but sufficient for our purposes. It is in the same style as the proof of \cite[Lemma 5.3]{Ancona}.
	
	\smallskip
	
	Let $A$ be an abelian variety and $\cL$ be an ample line bundle. There are two ways to associate with $\cL$ an alternating pairing on the Tate module $T(A).$ The first is by considering the Weil pairing between $T(A)$ and  the Tate module of the dual abelian variety $T(A^\vee)$ and then using $\cL$ to have an isogeny between $A$ and $A^\vee$ and hence identify $T(A)$ with $T(A^\vee)$. Notice that this pairing has the advantage of being non-degenerate (but it is unclear how to extend it to general commutative group schemes).
	
	The second pairing is that in the present paper. Look at the Chern class of $\cL$ in the $\ell$-adic cohomology $H^2(A)=\Lambda^2 H^1(A)$. By dualizing, such a class becomes an alternating   pairing on $T(A)$. 
	These two pairings are probably equal, but we do not know how to show it. 
	We instead show that they are the same up to an invertible scalar, which implies in particular that the second pairing is non-degenerate as well.
	
	To do so notice that, after fixing a suitable level structure, $A$ becomes a point of the moduli space  $M$  parameterizing abelian varieties of the same dimension as $A$ together with line bundles of the same degree as $\cL$. Now both pairings extend to the whole family of abelian varieties parameterized by $M$, hence define pairings at the level of the local systems. This means that both pairings on $T(A)$ are equivariant with respect to the monodromy action. Such an action on $T(A)$ is geometrically irreducible (it is the standard representation of the group $\mathrm{GSp}$); hence by Schur's lemma, the two pairings are the same up to a scalar. 
	Since both pairings are non-zero (the Chern class of an ample line bundle on a projective variety is non-zero), we conclude that such a scalar is invertible.
\end{proof}

The following is  \cite[Lemme XI.1.11]{Ray}. We warmly thank Michel Brion for pointing out  the original reference to us.
Nevertheless, we include a proof for completeness.

\begin{proposition}\label{p L e M ample}
	Let $\cL$ be an ample line bundle on $G$ which by \cref{p surj on Picard} is isomorphic to $p^*\cM$ with $\cM$ a line bundle on $A$. Then $\cM$ is ample. In other words, any ample line bundle on $G$ is the pullback of an ample line bundle on $A$.
\end{proposition}

\begin{proof}
	We know from \cite[Section~II.6, Application 1, p.~60]{MumfordTata} that $\cM$ is ample on $A$ if and only if the closed subgroup
	\[ K(\cM):=\{x\in A\mid t_x^*\cM\simeq \cM\} \]
	is finite.
	Suppose toward a contradiction that $K(\cM)$ is not finite, and let $Y$ be its identity component. Notice that $Y$ is itself an abelian variety. Denote by $\cM'=\cM|_Y$ the restriction of $\cM$ to $Y$.
	The proof from \textit{loc.\~cit.}~shows that on $Y$, we have 
	\begin{align}\label{eq -1 tensor triv}
		\cM'\otimes (-1_Y)^*\cM'\simeq \cO_Y.
	\end{align}
	
	On the subgroup $H:=p\inv Y$, we have that the ample line bundle $\cL|_H$ is isomorphic to $p^*\cM'$.
	However, if $p^*\cM'$ is ample, then  $(-1_H)^*p^*\cM'$ is also ample, and hence their tensor product is also ample.
	By \eqref{eq -1 tensor triv}, we get that $\cO_H$ is ample on $H$, \textit{i.e.}, that $H$ is quasi-affine, hence affine by \cref{l qaffine group is affine}.
	This gives a contradiction; therefore, $\cM$ is ample on $A$.	
\end{proof}


\newcommand{\etalchar}[1]{$^{#1}$}
\def\cprime{$'$}

\end{document}